\documentclass{amsart}[12pt]

\newtheorem{theorem}{Theorem}[section]
\newtheorem{lemma}[theorem]{Lemma}

\newtheorem{corollary}[theorem]{Corollary}

\theoremstyle{definition}
\newtheorem{definition}[theorem]{Definition}
\newtheorem{example}[theorem]{Example}

\usepackage{amscd,amssymb}

\begin{document}

\title{Free Assosymmetric algebras as modules of groups}
\author{Askar S. Dzhumadil'daev, Bekzat K. Zhakhayev}
\date{}
\address{Institute of Mathematics and Mathematical Modeling, Pushkin str., 125, Almaty, Kazakhstan}
\email{dzhuma@hotmail.com}

\address{Faculty of Engineering and Natural Sciences,
Suleyman Demirel University, Abylai Khan Str., 1/1, Kaskelen, 
Kazakhstan}
\email{bekzat22@hotmail.com}

\thanks
{The research of the first named author was supported by grant AP05131123 of the Ministry of Education and Science
of the Republic of Kazakhstan.}

\thanks
{The research of the second named author was supported by grant AP05131179 of the Ministry of Education and Science
of the Republic of Kazakhstan.}

\maketitle

\begin{abstract}
An algebra with identities $(a,b,c)=(a,c,b)=(b,a,c)$ is called {\it assosymmetric}, where $(x,y,z)=(xy)z-x(yz)$ is associator. We study $S_n$-module, $A_n$-module and $GL_n$-module structures of free assosymmetric algebra.
\end{abstract}

\section{Introduction}

Let $K$ be an algebraically closed field of characteristic $0$. All algebras, vector spaces, modules  and tensor products we consider will be over field $K$.  

Let $X=\{x_1, x_2,...\}$ be a set of generators and $K\{X\}$ be the  absolutely free nonassociative algebra.  A polynomial $f(x_1, x_2, ..., x_n)\in K\{X\}$ is called {\it polynomial identity} or {\it identity} for the $K$-algebra $R$ if $f(r_1, r_2,..., r_n)$=0 for all $r_1, r_2,..., r_n\in R$. 

Let $\{f_i\in K\{X\}| i\in I \}$ be a set of elements in $K\{X\}$. The class $\mathfrak{V}$ of all algebras satisfying the polynomial identities $f_i=0, i\in I$ is called the {\it variety} defined by the system of polynomial identities $\{f_i| i\in I\}$. The set  $T(\mathfrak{V})$ of all polynomial identities satisfied by the variety $\mathfrak{V}$ is called the $T$-{\it ideal} or {\it verbal ideal} of $\mathfrak{V}$.

A nonassociative algebra $({R}, \cdot)$ is called {\it  assosymmetric}, if for any
$a, b, c\in {R}$ the following identities are hold
\begin{displaymath}
(a, b, c)=(a, c, b)=(b, a, c),
\end{displaymath}
where $(x, y, z)=(x\cdot y)\cdot z-x\cdot(y\cdot z)$ is associator.

The factor algebra $F(X)=K\{X\}/ T({R})$ is called {\it free assosymmetric} algebra.

Assosymmetric algebras was studied in ~\cite{Klein},~\cite{PoRo},~\cite{Bo},~\cite{HentJaPe},~\cite{Dzhu}. 

Basis of free assosymmetric algebras was constricted in  ~\cite{HentJaPe}. 
Moreover, this paper contains multiplication rule of basis elements that allows to present an element of free assosymmetric algebra as a 
linear combination by basis elements. 
In ~\cite{Dzhu} it was proved that assosymmetric algebras under Jordan product satisfy Lie triple  and  Glennie identities. 


In polynomial identities theory there are two main  questions: $1)$ describe algebras with identities; $2)$ describe identities in algebras. The language of varieties allows one to freely pass from identity to algebra and from algebra to identity. Therefore studying varieties of algebras is one of the important problem in modern algebras. In 1950, A.I. Malcev ~\cite{Malcev} and W.Specht ~\cite{Specht} first time and independently used the representation theory of symmetric group to classify polynomial identities of algebraic structures. 
If $char K=0$, then every polynomial is equivalent to a finite set of multilinear polynomials. 

For several classes of algebras $S_n$-, $GL_n$-module structures on multilinear parts of free algebras are studied.  Some cases these structures can be easy described. 
For example, multilinear parts of free associative, free Zinbiel and free Leibniz algebras of degree $n$ as $S_n$-module are isomorphic to regular module $KS_n$~\cite{Zin}. 
In case of Lie algebras module structures are slightly complicated. 
In ~\cite{Klyach}  it was found list of irreducible $S_n$-representations that are involved in decomposition of multilinear parts of free Lie algebras.  
Description of multiplicities of irreducible $S_n$-representations in decomposition of multilinear part of free Lie algebra by language of major indices of standard Young tableaux is given in  ~\cite{Kras} .

Let us give some references for other papers where representations of  groups on free algebras are studied. 

In ~\cite{VlaDrens} full described varieties of associative algebras with identity of degree three by the methods of the theory of representations of symmetric group and general linear group. 
In ~\cite{AnaKem}  is given the criterion to the distributivity of the lattice of subvarieties of varieties of associative  algebras by using the methods of the theory of representations of symmetric group and general linear group.
In ~\cite{Mar} is given the criterion to the distributivity of the lattice of subvarieties of varieties of alternative algebras by using the methods of the theory of representations of symmetric group and general linear group.
In ~\cite{DrenZha} proved Specht problem (or finite basis problem) for varieties of bicommutative algebras over field of  characteristic $0$ by using methods of the theory of the representations of symmetric group.

In ~\cite{DzhuIsTu} constructed basis of free bicommutative algebras, described multiplication of basis elements, found cocharacter sequence, codimension sequence, calculated Hilbert series. It is also proved that bicommutative operad is not Koszul and the growth of codimension sequence of bicommutative algebra is equal to $2$. 
In ~\cite{Dren} is given alternative proof of the formula for the cocharacter sequence of bicommutative algebra.

In ~\cite{DzhuIs} studied module structures of free Novikov algebra over symmetric group and general linear group. It is given criteria for Novikov admissibility of any partition of positive integer.

In ~\cite{BremVar} M. Bremner  studied varieties of anticommutative $n$-ary algebras and found the "correct" generalization of the Jacobi identity. It is also formulated several conjectures. One of them is module structures of $P_n^2$ ($S_{2n-1}$-module spanned by the $n$-ary anticommutative monomials involving two pairs of brackets).
In ~\cite{Rot} M. Rotkiewicz proved above mentioned Bremner's conjecture.
In ~\cite{Brem} classified varieties of anti-commutative algebras defined by identities of degree $n\leq 7$. Author classified using computing the decomposition of the $S_n$-module of multilinear polynomials of degree n into irreducible submodules. Until today in science is unknown the module structures of free anti-commutative algebra in general case. 

In ~\cite{Labelle} calculated characters of representations of symmetric group on free right-symmetric and right-commutative algebras respectively.

In ~\cite{HenRegExp} showed explicitly decomposition of the group algebra of the alternating group into direct sum of minimal left ideals.
In ~\cite{HenRegWeyl} studied Weyl modules for the Schur algebra of the Alternating Group.
In ~\cite{HenReg} calculated the $A$-codimensions and the $A$-cocharacters of the infinite dimensional Grassmann (exterior) algebra. Authors conjectured a finite generating set of the $A_n$-identities for the Grassmann algebra. 
In ~\cite{GonKosh} proved Henke-Regev conjecture.

In this paper we consider multilinear part of free
assosymmetric algebra $F(X)$ as $S_n$-module, $A_n$-module and $GL_n$-module. We find dimension of homogeneous component, sequence of  dimensions of multilinear components or codimension sequence, colength sequence, cocharacter sequence in $S_n$-case, cocharacter sequence in $A_n$-case for assosymmetric algebras.

\section{Statement of main result}

Let $n$ be a positive integer. The sequence of positive integers $\lambda=(\lambda_1, \lambda_2, ..., \lambda_k)$ is called {\it partition} of $n$, if (1) $\lambda_1+\lambda_2+...+\lambda_k=n$, (2) $\lambda_1\geq \lambda_2\geq ... \geq \lambda_k$, and denoted by $\lambda\vdash n$. {\it Length} of partition $\lambda\vdash n$ is the number of parts in $\lambda$ and denoted by $\ell(\lambda)$. It is known that between partitions of $n$ and Young diagrams with $n$ boxes exist one-to-one correspondence. We denote Young diagram with $\lambda$-shape by $Y_{\lambda}$. Let $\lambda, \mu\vdash n$. Partition $\lambda$ is {\it conjugate} to partition $\mu$, if $Y_{\mu} $ is obtained from $Y_{\lambda}$ by turning the rows into columns and denoted by $\mu=\lambda'$.  A partition that is conjugate to itself is said to be a {\it self-conjugate} partition, that is $\lambda=\lambda'$. 

Let $S_n$ be symmetric group on set $\{1,2,...,n\}$ and $A_n$ be alternating subgroup of $S_n$. The symmetric group $S_n$ and alternating group $A_n$ acts on multilinear part of free assosymmetric algebra in natural way (left action or variable action). 

Let $R$ be a PI-algebra. For $n\geq 1$, the $S_n$-character of $P_n(R)=P_n/(P_n\cap T(R))$ is called the {\it $n$-th cocharacter} of $R$ and denoted by $\chi_n(R)$, and 
\begin{displaymath}
\chi_n (R)=\sum_{\lambda \vdash n} m_{\lambda} \chi_{\lambda},
\end{displaymath}
where $\chi_{\lambda}$ is the irreducible $S_n$-character associated to the partition $\lambda\vdash n$ and $m_{\lambda}\geq 0$ is the corresponding multiplicity.

Let $T(R)$ be $T$-idel of $R$. Then the non-negative integer
\begin{displaymath}
c_n(R)=dim (P_n/P_n\cap T(R)),
\end{displaymath}
is called the {\it $n$-th codimension} of the algebra $R$.

Let $R$ be a PI-algebra and 
\begin{displaymath} 
\chi_n(R)=\sum_{\lambda\vdash n} m_{\lambda} \chi_{\lambda}.
\end{displaymath}
Then the non-negative integer
\begin{displaymath}
l_n(R)=\sum_{\lambda\vdash n} m_{\lambda}
\end{displaymath}
is called the {\it $n-$th colength} of $R$.

For more information about codimension sequence, cocharacter sequence, colength sequence see ~\cite{GiamZai}.

We denote irreducible $S_n$-module or Specht module associated to partition $\lambda\vdash n$ by $S^{\lambda}$ and dimension of $S^{\lambda}$ by $d_{\lambda}$, irreducible $A_n$-module associated to non-self-conjugate partition $\lambda\vdash n$ by $S^{\lambda}_A$ and to self-conjugate partition $\lambda\vdash n$ by $S^{\lambda\pm}_A$, irreducible $GL_n$-module or Weyl module associated to partition $\lambda\vdash n$ by $W^{\lambda}$ in $S_n$-case and by $W^{\lambda}_A$ in $A_n$-case. 

For more information about the theory of representations of $S_n$, $A_n$ and $GL_n$ see ~\cite{Ful}, ~\cite{FulHar} ~\cite{JamesKer}, ~\cite{Sagan}, ~\cite{HenReg}, ~\cite{HenRegExp}, ~\cite{HenRegWeyl}

Free base of assosymmetric algebras was found  in \cite{HentJaPe}. We use this result to find formulas for dimensions of free assosymmetric algebras. 
Let $F(r)$ be free assosymmetric algebra generated by $r$ elements $a_1,\ldots,a_r.$ 
Let  $F^{l_1,\ldots,l_r}(r)$ be a subspace of free assosymmetric algebra generated by $l_i$ elements $a_i,$ where $i=1,\ldots,r,$ and 
$F_n(r)$ be a subspace of free assosymmetric algebra $F(r)$ of degree $n$  and $F_n^{multi}=F^{1,\ldots, 1}(n)$ be multi-linear part of $F_n(n).$

\begin{theorem} \label{dim} Let $p=char K\ne 2,3.$ Then 
$$dim\,F^{l_1,\ldots,l_r}(r)=
{l_1+\cdots+l_r\choose l_1\; \cdots \; l_r}+(l_1+1)\cdots(l_r+1)-{r+1\choose 2}-r-1+w,$$
where $w=w(l_1,\ldots,l_r)$ is a number of $1$'s in the sequence $l_1\ldots l_r,$

\medskip

$$dim\,F_n(r)=$$ 
$$r^n+{n+2r-1\choose n}-{r+1\choose 2}{n+r-3\choose n-2}-r {n+r-2\choose n-1}-{n+r-1\choose n},$$

\medskip

\noindent and

\medskip

$$dim\,F_n^{multi}=n!+2^n-{n+1\choose 2}-1.$$
$\square$
\end{theorem}

By Stirling formula $n!\sim  \sqrt{2\pi n}(n/e)^n,$ and therefore,
$$ dim\,F_n^{1/n} \sim n/e .$$

We divide the set of multilinear basic elements into two types.

First type
\begin{displaymath}
T_n=\{(...((x_{\sigma(1)} x_{\sigma(2)})x_{\sigma(3)})...)x_{\sigma(n)}| \sigma\in S_n\},
\end{displaymath}

Second type
\begin{displaymath}
T_{k, n-k}=\{x_{\sigma(1)} (x_{\sigma(2)}( ... x_{\sigma(k)}[... [(x_{\sigma(k+1)} x_{\sigma(k+2)}, x_{\sigma(k+3)}), x_{\sigma(k+4)}],..., x_{\sigma(n)}]...))| 
\end{displaymath}
\begin{displaymath}
\sigma(1) < \sigma(2) <  ... < \sigma(k), \  \ \sigma(k+1)< \sigma(k+2)< ... < \sigma(n)\}.
\end{displaymath}

\begin{theorem}\label{thm2}
The group $S_n$ acts transitively on the sets $T_n$ and $T_{k, n-k}$, \ \ $k=0,1,...,n-3$.
\end{theorem}

Let $KT_n$ and  $KT_{k,n-k}$ be subspaces of $F_n^{multi}:=P_n$ spanned by the sets $T_n$ and $T_{k,n-k}$, $k=0,1,...,n-3$, respectively.

\begin{corollary}\label{cor1}
As $S_n$-module
\begin{displaymath}
P_n\cong KT_n\oplus \bigoplus_{k=0,1,...,n-3} KT_{k, n-k}.
\end{displaymath}
\end{corollary}

\begin{theorem}\label{thm3}
As $S_n$-module
\begin{displaymath}
P_n \cong \bigoplus_{\lambda \vdash n} d_{\lambda} S^{\lambda} \oplus \bigoplus_{(\lambda_1, \lambda_2)\vdash n} m(\lambda_1, \lambda_2) S^{(\lambda_1, \lambda_2)},
\end{displaymath}
where
\begin{displaymath}
m(\lambda_1, \lambda_2)=
\left\{\begin{array}{ll}
n-2-\lambda_2, & \lambda_2\leq 3,\\
n+1-2\lambda_2, & \lambda_2\geq 4.
\end{array}
\right.
\end{displaymath}
\end{theorem}

\begin{example}
\begin{displaymath}
P_1\cong S^{(1)};
\end{displaymath}
\begin{displaymath}
P_2\cong S^{(2)}\oplus S^{(1,1)};
\end{displaymath}
\begin{displaymath}
P_3\cong 2*S^{(3)}\oplus 2*S^{(2,1)}\oplus S^{(1,1,1)};
\end{displaymath}
\begin{displaymath}
P_4\cong 3*S^{(4)}\oplus 4*S^{(3,1)}\oplus 2*S^{(2,2)}\oplus 3*S^{(2,1,1)}\oplus S^{(1,1,1,1)};
\end{displaymath}
\begin{displaymath}
P_5\cong 4*S^{(5)}\oplus 6*S^{(4,1)}\oplus 6*S^{(3,2)}\oplus 6*S^{(3,1,1)}\oplus 5*S^{(2,2,1,1)}\oplus 4*S^{(2,1,1,1,1)}\oplus S^{(1,1,1,1,1)}.
\end{displaymath}
\end{example}

Let $V$ be a vector space with dimension $m$. Let $F(V)$ be free assosymmetric algebra generated by basis elements of $V$ and $H_n(V)$ be homogeneous part of $F(V)$ of degree $n$.

\begin{definition} Let $inv(S_n)$ be number of involutions, 
\begin{displaymath}
inv(S_n)= \# \{ \sigma\in S_n \ \ |\ \  \sigma^2=e\}
\end{displaymath}
\end{definition}

\begin{corollary}\label{cor2}
{\bf a.} \begin{displaymath}
\chi_{S_n}(P_n)\cong \sum_{\lambda \vdash n}d_{\lambda} \chi_{S_n}(\lambda)+\sum_{(\lambda_1, \lambda_2)\vdash n} m(\lambda_1, \lambda_2) \chi_{S_n}(\lambda_1, \lambda_2),
\end{displaymath}
where 
\begin{displaymath}
m(\lambda_1, \lambda_2)=
\left\{\begin{array}{ll}
n-2-\lambda_2, & \lambda_2\leq 3,\\
n+1-2\lambda_2, & \lambda_2\geq 4.
\end{array}
\right.
\end{displaymath}

{\bf b.}
\begin{displaymath}
\chi_{A_n}(P_n)=2\chi_{A_n}(KA_n)+\sum_{(\lambda_1, \lambda_2)\vdash n} m(\lambda_1, \lambda_2) \chi_{A_n}(\lambda_1, \lambda_2), 
\end{displaymath}
where 
\begin{displaymath}
m(\lambda_1, \lambda_2)=
\left\{\begin{array}{ll}
n-2-\lambda_2, & \lambda_2\leq 3,\\
n+1-2\lambda_2, & \lambda_2\geq 4.
\end{array}
\right.
\end{displaymath}

{\bf c.} ($S_n$-case)
\begin{displaymath}
H_n(V)\cong \bigoplus_{\lambda\vdash n} d_{\lambda} W^{\lambda}\oplus \bigoplus_{(\lambda_1, \lambda_2)\vdash n} m(\lambda_1, \lambda_2) W^{\lambda},
\end{displaymath}
where 
\begin{displaymath}
d_{\lambda}>0,\ \ m(\lambda_1, \lambda_2)>0 \ \ \ \ \text{and}
\end{displaymath}
\begin{displaymath}
m(\lambda_1, \lambda_2)=
\left\{\begin{array}{ll}
n-2-\lambda_2, & \lambda_2\leq 3,\\
n+1-2\lambda_2, & \lambda_2\geq 4,
\end{array}
\right.
\end{displaymath}
if $dim V\geq \ell (\lambda), \ell ((\lambda_1, \lambda_2))$, and

\begin{displaymath}
d_{\lambda}=0, \ \ m(\lambda_1, \lambda_2)=0,
\end{displaymath}
if $ dim V <  \ell (\lambda), \ell((\lambda_1, \lambda_2))$.

{\bf d.} ($A_n$-case)
\begin{displaymath}
H_n(V)\cong \left[\bigoplus_{\lambda\neq \lambda'}2 d_{\lambda}W^{\lambda}_A \right] \oplus \left[\bigoplus_{\lambda=\lambda'}2\left(\frac{d_{\lambda}}{2}W^{\lambda+}_A\oplus \frac{d_{\lambda}}{2} W^{\lambda-}_A\right)\right]\oplus \bigoplus_{(\lambda_1, \lambda_2)\vdash n} m(\lambda_1, \lambda_2) W^{\lambda}_A,
\end{displaymath}
where 
\begin{displaymath}
m(\lambda_1, \lambda_2)=
\left\{\begin{array}{ll}
n-2-\lambda_2, & \lambda_2\leq 3,\\
n+1-2\lambda_2, & \lambda_2\geq 4.
\end{array}
\right.
\end{displaymath}

{\bf e.} For $1\leq n\leq3$
\begin{displaymath}
l_n(P_n)=\delta_{n,3}+inv(S_n),
\end{displaymath}
where $\delta_{i,j}$ is Kronecker delta.

For $n\geq 4$
\begin{displaymath}
l_n(P_n)=
\left\{\begin{array}{ll}
k^2+2k-5+inv(S_n), & n=2k,\\
k^2+3k-4+inv(S_n), & n=2k+1.
\end{array}
\right.
\end{displaymath}
\end{corollary}

\begin{example}
\begin{displaymath}
\chi_{S_1}(P_1)=\chi_{S_1}{(1)};
\end{displaymath}
\begin{displaymath}
\chi_{S_2}(P_2)=\chi_{S_2}{(2)}+ \chi_{S_2}{(1,1)};
\end{displaymath}
\begin{displaymath}
\chi_{S_3}(P_3)=2*\chi_{S_3}{(3)}+ 2*\chi_{S_3}{(2,1)}+ \chi_{S_3}{(1,1,1)};
\end{displaymath}
\begin{displaymath}
\chi_{S_4}(P_4)=3*\chi_{S_4}{(4)}+ 4*\chi_{S_4}{(3,1)}+ 2*\chi_{S_4}{(2,2)}+ 3*\chi_{S_4}{(2,1,1)}+ \chi_{S_4}{(1,1,1,1)};
\end{displaymath}
\begin{displaymath}
\chi_{S_5}(P_5)=4*\chi_{S_5}{(5)}+ 6*\chi_{S_5}{(4,1)}+ 6*\chi_{S_5}{(3,2)}+ 6*\chi_{S_5}{(3,1,1)}+ 5*\chi_{S_5}{(2,2,1)}+
\end{displaymath}
\begin{displaymath}
+4*\chi_{S_5}{(2,1,1,1,1)}+ \chi_{S_5}{(1,1,1,1,1)}.
\end{displaymath}
\end{example}

\begin{example}
\begin{displaymath}
\chi_{A_1}(P_1)=\chi_{A_1}{(1)};
\end{displaymath}
\begin{displaymath}
\chi_{A_2}(P_2)=\chi_{A_2}{(2)}+ \chi_{A_2}{(1,1)};
\end{displaymath}
\begin{displaymath}
\chi_{A_3}(P_3)=3*\chi_{A_3}{(3)}+ \chi_{A_3}^{+}{(2,1)}+\chi_{A_3}^{-}{(2,1)};
\end{displaymath}
\begin{displaymath}
\chi_{A_4}(P_4)=4*\chi_{A_4}{(4)}+ 7*\chi_{A_4}{(3,1)}+\chi_{A_4}^{+}{(2,2)}+\chi_{A_4}^{-}{(2,2)};
\end{displaymath}
\begin{displaymath}
\chi_{A_5}(P_5)=5*\chi_{A_5}{(5)}+ 10*\chi_{A_5}{(4,1)}+ 11*\chi_{A_5}{(3,2)}+ 3*\chi_{A_5}^{+}{(3,1,1)}+3*\chi_{A_5}^{-}{(3,1,1)}.
\end{displaymath}
\end{example}

\begin{example}
\begin{displaymath}
l_1(P_1)=1; \ \ l_2(P_2)=2; \ \ l_3(P_3)=5; \ \ l_4(P_4)=13;\ \ l_5(P_5)=32.
\end{displaymath}
\end{example}

\section{Proof of Theorem \ref{dim}}

In calculation of dimensions we need the following easy proved combinatorial results. 

\begin{lemma} \label{ppp}
For non-negative integers $\alpha,\beta$ and $n$ takes place the following formula
$$\sum_{i=0}^n {i+\alpha\choose i}{n-i+\beta\choose n-i}={n+\alpha+\beta+1\choose n}$$
In particular,
$$\sum_{i=0}^n {i+\alpha\choose i}{n-i+\alpha\choose n-i}={n+2\alpha+1\choose n}$$
\end{lemma}

\begin{lemma}\label{pp}
 Number of non-decreasing sequences of length $m$ with components in the set $I=\{1,2,\ldots,r\}$ is ${m+r-1\choose m}.$
\end{lemma}

\begin{lemma}\label{qqqq}
 Number of non-decreasing sequences of length $m$ with components in the set $I=\{1,2,\ldots,r\}$  such that each $i\in I$  appears no more than $l_i$ times is $(l_1+1)\cdots (l_r+1).$
 \end{lemma}

In \cite{HentJaPe} is proved that a base of free assosymmetric algebras can be constructed by elements of two kinds.  
 If $S=\{a_1,\ldots,a_r\}$ is a set of  generators, then in degree $n$ the base consists elements of a form 
$$(\cdots ((a_{i_1}a_{i_2})a_{i_3})\cdots) a_{i_n}, \quad a_i\in S,$$
$$a_{i_1}(a_{i_2}(\cdots a_{i_m}[\cdots[(a_{j_1},a_{j_2},a_{j_3}),a_{j_4}],\cdots,a_{j_k}]\cdots)), \quad a_i, a_j\in S,$$
$$i_1\le i_2\le \cdots\le i_m, \;\;  j_1\le j_2\le \cdots \le j_k, \quad  m\ge 0, \; k\ge 3.$$
Number of elements of first kind is $r^n.$ By Lemma \ref{pp} number of elements of second kind $L$ is equal to 
$$L=\sum_{m+k=n, m\ge 0, k\ge 3} {m+r-1\choose m}{k+r-1\choose k}=$$
$$\sum_{m+k=n, m\ge 0, k\ge 0}{m+r-1\choose m}{k+r-1\choose k}-$$
$${n+r-3\choose n-2}{r+1\choose 2}-{n+r-2\choose n-1}{r\choose 1}-{n+r-1\choose n}{r-1\choose 0}$$
By Lemma \ref{ppp}
$$L={n+2r-1\choose n}-{r+1\choose 2}{n+r-3\choose n-2}-r{n+r-2\choose n-1}-{n+r-1\choose n}.$$
Therefore, 
$$dim\,F_n(r)=$$ $$r^n+
{n+2r-1\choose n}-{r+1\choose 2}{n+r-3\choose n-2}-r{n+r-2\choose n-1}-{n+r-1\choose n}.$$

Now suppose that any generator $a_s,$ $s=1,2,\ldots,r,$ in each base element should enter $l_s$ times.  Then the  number of base elements of first kind is $${l_1+\cdots+l_n \choose l_1 \, \cdots \, l_n}=\frac{(l_1+\cdots + l_r)!}{l_1!\cdots l_r!}.$$ 

Let  $M$ be set of sequences $\alpha=i_1  \ldots i_m j_1 j_2\ldots j_k$  with components in $I=\{1,2,\ldots,r\}$  such that each $i\in I$ appears exactly $l_i$ times and 
$i_1\le \cdots \le i_m,$ $j_1\le \cdots \le j_k.$ 
For $\alpha\in M$ call its subsequence of first $m$ components $i_1\ldots i_m$ as {\it head} 
and denote $\tilde \alpha.$
Note that each $\alpha\in M$ 
is uniquely defined by head $\tilde \alpha.$ 
Denote set of heads  by $\tilde M.$
Note also that in the sequence $\tilde\alpha=i_1\ldots i_m$ each $i\in I$ enters no more than  $l_i$ times. Therefore by Lemma \ref{qqqq} the number of heads is  
$$|\tilde M|=(l_1+1)\cdots (l_r+1).$$

Let $N$ be a subset of $M$  consisting of sequences with the following heads 
$$\underbrace{1\ldots 1}_{l_1}\ldots \underbrace{i\ldots i}_{l_i} \ldots \underbrace{r\ldots r}_{l_r},$$
(number of  such sequences is $1$)

$$\underbrace{1\ldots 1}_{l_1}\ldots \underbrace{i\ldots i}_{l_i-1} \ldots \underbrace{r\ldots r}_{l_r},\quad i\in I,$$
(number of  such sequences is $r$)

$$\underbrace{1\ldots 1}_{l_1}\ldots \underbrace{i\ldots i}_{l_i-2} \ldots \underbrace{r\ldots r}_{l_r}, \quad l_i>1, \quad i\in I,$$
(number of  such sequences is $r-w$, where $w$ is a number of $1$'s in the sequence 
$l_1\ldots l_r$)

$$\underbrace{1\ldots 1}_{l_1}\ldots
\underbrace{i\ldots i}_{l_i-1} \ldots
\underbrace{j\ldots j}_{l_j-1}\ldots
 \underbrace{r\ldots r}_{l_r}, \quad i<j, \quad i,j\in I$$
(number of  such sequences is $r(r-1)/2$).

Let $M_1=M\setminus N$ be a supplement of $N$ in the set $M.$ Then any $\alpha= i_1\ldots i_m j_1\ldots j_k \in M_1$ has the property $k\ge 3$ and any such sequence generates base element of free assosymmetric algebra of second kind. Hence the  number of base elements of second kind is 
$$dim\,F^{l_1,\ldots,l_r}(r)=$$
$$|M_1|=
{l_1+\cdots+l_r\choose l_1\; \cdots \; l_r}+(l_1+1)\cdots(l_r+1)-{r+1\choose 2}-r-1+w.$$
Dimension for multilinear part is an easy consequence of this formula.

\section{Proof of Theorem \ref{thm2}}

Let $A$ be the $T$-ideal in $K\{X\}$ determined by identities

\begin{displaymath}
(x,y,z)=(x,z,y)=(y,x,z).
\end{displaymath}

Let $R$ be assosymmetric algebra and $I=(R,R,R)+(R,R,R)R$ is the ideal generated by associators.
Proof of Theorem \ref{thm2} is based on the following two results of \cite{HentJaPe}. 

\begin{lemma}[~\cite{HentJaPe}, Lemma 1] \label{Bekzat1}
The expression $[[...[[(a_1,a_2,a_3),a_4],a_5]...]a_n]$ is invariant, modulo A, under all permutations of the arguments.
\end{lemma}

\begin{lemma}[~\cite{HentJaPe}, Lemma 2]\label{Bekzat2}
If $x\in I$, the expression $a_1(a_2(a_3(... a_n x)...))$ is invariant, modulo A, under all permutations of the $a_i$'s.
\end{lemma}

Now we give proof of Theorem \ref{thm2}.

First type. Let
\begin{displaymath}
(...((a_{i_1}\cdot a_{i_2})\cdot a_{i_3})\cdot...)\cdot a_{i_n}\in T_n, \  \ i_j\in \{1,2,...,n\}.
\end{displaymath}
Let $\sigma=(i_k, i_{k+1})$ be a transposition in $S_n$. Then 
\begin{displaymath}
(i_k, i_{k+1}): (...((a_{i_1}\cdot a_{i_2})\cdot a_{i_3})\cdot...)\cdot
a_{i_k})\cdot a_{i_{k+1}})\cdot...)\cdot a_{i_n}\mapsto
(...((a_{i_1}\cdot a_{i_2})\cdot a_{i_3})\cdot...)\cdot
a_{i_{k+1}})\cdot a_{i_{k}})\cdot...)\cdot a_{i_n}.
\end{displaymath}
By definition of  $(...((a_{i_1}\cdot a_{i_2})\cdot
a_{i_3})\cdot...)\cdot a_{i_{k+1}})\cdot a_{i_{k}})\cdot...)\cdot
a_{i_n}$ is multilinear basis element in $T_n$.

Second type. Let
\begin{displaymath}
w= a_{i_1}\cdot (a_{i_2}\cdot( ... a_{i_k}\cdot[... [(a_{i_{k+1}},
a_{i_{k+2}}, a_{i_{k+3}}), a_{i_{k+4}}],..., a_{i_n}]...))\in T_{k, n-k}, \ \ i_j\in \{1,2,...,n\}.
\end{displaymath}
We present it in a form 

\begin{displaymath}
\underbrace{a_{i_1}\cdot (a_{i_2}\cdot( ...
a_{i_k}}_{A-part}\cdot\underbrace{[... [(a_{i_{k+1}}, a_{i_{k+2}},
a_{i_{k+3}}), a_{i_{k+4}}],..., a_{i_n}]}_{B-part}...)).
\end{displaymath}


It suffices to consider the action of transposition $\sigma=(i_j, i_{j+1})\in S_n$ in three cases:

{\it Case-1:} $\sigma$ acts on $A$-part;

{\it Case-2:} $\sigma$ acts on $B$-part;

{\it Case-3:} $\sigma$ acts on $A$-part and $B$-part simultaneously.

{\it Case-1:}

\begin{displaymath}
\sigma: a_{i_1} (a_{i_2}( ... a_{i_j}(a_{i_{j+1}} ... (
a_{i_k} [... [(a_{i_{k+1}}, a_{i_{k+2}},
a_{i_{k+3}}), a_{i_{k+4}}],..., a_{i_n}]...))\mapsto 
\end{displaymath}
\begin{displaymath}
\mapsto a_{i_1} (a_{i_2}( ... a_{i_{j+1}}(a_{i_j} ... (
a_{i_k} [... [(a_{i_{k+1}}, a_{i_{k+2}},
a_{i_{k+3}}), a_{i_{k+4}}],..., a_{i_n}]...)).
\end{displaymath}

By Lemma \ref{Bekzat2} \ \ $\sigma w\in T_{k, n-k}$  and  $\sigma w=w$.

{\it Case-2:}

\begin{displaymath}
\sigma: a_{i_1} (a_{i_2}( ... (
a_{i_k} [... [(a_{i_{k+1}}, a_{i_{k+2}},
a_{i_{k+3}}), a_{i_{k+4}}],...,a_j], a_{j+1},..., a_{i_n}]...))\mapsto 
\end{displaymath}
\begin{displaymath}
\mapsto a_{i_1} (a_{i_2}( ... (
a_{i_k} [... [(a_{i_{k+1}}, a_{i_{k+2}},
a_{i_{k+3}}), a_{i_{k+4}}],...,a_{j+1}], a_j,...,  a_{i_n}]...)).
\end{displaymath}

By Lemma \ref{Bekzat1} \ \ $\sigma w\in T_{k, n-k}$  and  $\sigma w=w$.

{\it Case-3:}

Let

\begin{displaymath}
w= a_{i_1}\cdot (a_{i_2}\cdot( ... a_{i_k}\cdot[... [(a_{i_{k+1}},
a_{i_{k+2}}, a_{i_{k+3}}), a_{i_{k+4}}],..., a_{i_n}]...)) \in T_{k, n-k}, \ \ i_j\in \{1,2,...,n\}.
\end{displaymath}

Assume that $a_{i_j}$ belong to $A$-part and $a_{i_{j+1}}$ belong to $B$-part, i.e.

\begin{displaymath}
w= a_{i_1} (a_{i_2} ( ...( a_{i_j}(...( a_{i_k}[... [(a_{i_{k+1}}, a_{i_{k+2}}, a_{i_{k+3}}), a_{i_{k+4}}],...]a_{i_{j+1}}]..., a_{i_n}]...)).
\end{displaymath}
{Then}

\begin{displaymath}
\sigma: a_{i_1} (a_{i_2} ( ...(a_{i_{j}}(...( a_{i_k} [... [(a_{i_{k+1}}, a_{i_{k+2}}, a_{i_{k+3}}), a_{i_{k+4}}],...]a_{i_{j+1}}]..., a_{i_n}]...))\mapsto
\end{displaymath}
\begin{displaymath}
\mapsto a_{i_1} (a_{i_2} ( ...(a_{i_{j+1}}(...( a_{i_k} [... [(a_{i_{k+1}}, a_{i_{k+2}}, a_{i_{k+3}}), a_{i_{k+4}}],...]a_{i_{j}}]..., a_{i_n}]...))
\end{displaymath}
\begin{displaymath}
=(\text{by Lemma \ref{Bekzat2} and Lemma \ref{Bekzat1}})
\end{displaymath}
\begin{displaymath}
=a_{p_1} (a_{p_2} ( ...(a_{p_{l}}(...( a_{p_k} [... [(a_{p_{k+1}}, a_{p_{k+2}}, a_{p_{k+3}}), a_{p_{k+4}}],...]a_{p_{m}}]..., a_{p_n}]...)),
\end{displaymath}
 where $\{i_1, i_2,\dots , i_{j+1},\dots , i_k\}=\{p_1, p_2,\dots , p_l,\dots , p_k\ \ | \ \ p_1<p_2<...<p_l<...<p_k\}$ and $\{i_{k+1}, i_{k+2}, i_{k+3}, i_{k+4}\dots , i_j, \dots, i_n\}=\{p_{k+1}, p_{k+2}, p_{k+3}, ... ,p_{m},...p_{n}\ \ | \ \ p_{k+1}<p_{k+2}<p_{k+3}<...<p_{m}<...<p_{n}\}$.

As we have noticed $\sigma w\neq w$.\ \ $\Box$

\section{Proof of Theorem \ref{thm3}}
Let $V(n)$ be a vector space with dimension $n$. By Theorem \ref{thm2} $KT_n$ is isomorphic to $\underbrace{V(1)\otimes V(1)\otimes ... \otimes V(1)}_{n}$ as $S_n$-module. Therefore

\begin{displaymath}
KT_n \cong Ind^{S_n}_{S_1\times S_1\times\dots \times S_1}({\bf 1}_{S_1}\otimes {\bf 1}_{S_1}\otimes\dots \otimes {\bf 1}_{S_1}) \cong \bigoplus_{\lambda\vdash n} d_{\lambda} S^{\lambda}, 
\end{displaymath} 
where ${\bf 1}_{S_1}$ is one-dimensional trivial representation of $S_1$.

{By Theorem \ref{thm2}}  group of automorphisms of $A$-part of $T_{k, n-k}$ is $S_k$ and group of automorphisms of $B$-part is $S_{n-k}$. Therefore $S_k\times S_{n-k}$ is group of automorphisms of $T_{k, n-k}$. 

Let
\begin{displaymath}
g_A=\sum_{\sigma\in S_k}\sigma\in KS_k, \ \ 
g_B=\sum_{\tau\in S_{n-k}}\tau\in KS_{n-k}.
\end{displaymath}
be elements of group algebras $KS_k$ and $KS_{n-k}$, respectivley. Then {by Theorem \ref{thm2}} $g_{T_{k, n-k}}=g_A\otimes g_B$ is  generator of {all} basic elements of $KT_{k, n-k}$ and $g_A,  g_B$ are {one-dimensional} trivial representations of $S_k$ and $S_{n-k}$, respectively, and $KT_{k, n-k}$ is $S_k\times S_{n-k}$-module. Therefore $KT_{k, n-k}$ as $S_n$-module is isomorphic to

\begin{displaymath}
Ind^{S_n}_{S_k\times S_{n-k}}({\bf 1}_{S_k}\otimes {\bf 1}_{S_{n-k}})\cong \bigoplus_{(\lambda_1, \lambda_2)\vdash n} S^{(\lambda_1, \lambda_2)}, \ \ \lambda_2\leq min\{k, n-k\},
\end{displaymath}
where ${\bf 1}_{S_k}=g_A$, ${\bf 1}_{S_{n-k}}=g_B$.

By Corollary \ref{cor1}

\begin{displaymath}
P_n\cong KS_n\oplus \bigoplus_{k=0,1,...,n-3} KT_{k, n-k}\cong
\end{displaymath} 
\begin{displaymath}
\cong \bigoplus_{\lambda \vdash n} d_{\lambda} S^{\lambda} \oplus \bigoplus_{(\lambda_1, \lambda_2)\vdash n} m(\lambda_1, \lambda_2) S^{(\lambda_1, \lambda_2)},
\end{displaymath}
where
\begin{displaymath}
m(\lambda_1, \lambda_2)=
\left\{\begin{array}{ll}
n-2-\lambda_2, & \lambda_2\leq 3,\\
n+1-2\lambda_2, & \lambda_2\geq 4.
\end{array}
\right.  \ \ \ \ \Box
\end{displaymath}

\section{Proof of Corollary \ref{cor2} }

 {\bf a.} Follows from Theorem \ref{thm3}. 

{\bf b.} $KT_n$ as $S_n$-module is isomorphic to

\begin{displaymath}
KT_n \cong \bigoplus_{\lambda\vdash n}d_{\lambda}S^{\lambda}.
\end{displaymath}

$KA_n$ as $A_n$-module is isomorphic to

\begin{displaymath}
KA_n\cong \left[\bigoplus_{\lambda\neq \lambda'}d_{\lambda} S^{\lambda}_A\right]\oplus \left[\bigoplus_{\lambda=\lambda'}(\frac{d_{\lambda}}{2}S^{\lambda+}_A\oplus \frac{d_{\lambda}}{2} S^{\lambda-}_A)\right],
\end{displaymath}
where $S_A^{\lambda}$ is irreducible $A_n$-module.

If $\lambda\vdash n$ is non-self-conjugate partition, then $S^{\lambda}$ and $S^{\lambda'}$ as $A_n$-modules are isomorphic to

\begin{displaymath}
Res^{S_n}_{A_n}(S^{\lambda})\cong S^{\lambda}_{A}, \ \ \ Res^{S_n}_{A_n}(S^{\lambda'})\cong S^{\lambda'}_{A}
\end{displaymath}
and
\begin{displaymath}
S^{\lambda}_A\cong S^{\lambda'}_A,
\end{displaymath}
where $dim(S^{\lambda}_A)=dim(S^{\lambda'}_A)=d_{\lambda}$.

If $\lambda\vdash n$ is self-conjugate partition, then $S^{\lambda}$ as $A_n$-module is isomorphic to

\begin{displaymath}
Res^{S_n}_{A_n}(S^{\lambda})\cong (S^{\lambda+}_{A}\oplus S^{\lambda-}_{A}),
\end{displaymath}
where $dim(S^{\lambda+}_A)=dim(S^{\lambda-}_A)=\frac{d_{\lambda}}{2}$. {Details } see ~\cite{HenRegExp}.

Therefore

\begin{displaymath}
KT_n \cong 2\cdot KA_n.
\end{displaymath}

{Note that} $KT_{k,n-k}$, $k=0,1,...,n-3$, as $S_n$-module is isomorphic to

\begin{displaymath}
KT_{k,n-k}\cong \bigoplus_{(\lambda_1, \lambda_2)\vdash n} S^{(\lambda_1, \lambda_2)}, \ \ \ \ \lambda_2\leq min\{k, n-k\}.
\end{displaymath}

Therefore $KT_{k,n-k}$ as $A_n$-module is isomorphic to

\begin{displaymath}
Res^{S_n}_{A_n}(KT_{k,n-k})\cong Res^{S_n}_{A_n}(\bigoplus_{(\lambda_1, \lambda_2)\vdash n} S^{(\lambda_1, \lambda_2)})\cong 
\end{displaymath}
\begin{displaymath}
\cong \bigoplus_{(\lambda_1, \lambda_2)\vdash n}Res_{A_n}^{S_n} S^{(\lambda_1, \lambda_2)}\cong \bigoplus_{(\lambda_1, \lambda_2)\vdash n} S^{(\lambda_1, \lambda_2)}_{A}, \ \ \lambda_2\leq min\{k, n-k\}.
\end{displaymath}

{\bf c.} ($S_n$-case) It is well known, that

\begin{displaymath}
W^{\lambda}\cong V^{\otimes n}\otimes_{KS_n}S^{\lambda}.
\end{displaymath}
Then
\begin{displaymath}
H_n(V)\cong V^{\otimes n}\otimes_{KS_n}P_n\cong
\end{displaymath}
\begin{displaymath}
\cong V^{\otimes n}\otimes_{KS_n}\left( \bigoplus_{\lambda \vdash n} d_{\lambda} S^{\lambda} \oplus \bigoplus_{(\lambda_1, \lambda_2)\vdash n} m(\lambda_1, \lambda_2) S^{(\lambda_1, \lambda_2)}\right)\cong
\end{displaymath}
\begin{displaymath}
\cong \left(V^{\otimes n}\otimes_{KS_n} \bigoplus_{\lambda \vdash n} d_{\lambda} S^{\lambda}\right) \oplus \left( V^{\otimes n}\otimes_{KS_n} \bigoplus_{(\lambda_1, \lambda_2)\vdash n} m(\lambda_1, \lambda_2) S^{(\lambda_1, \lambda_2)}\right)\cong
\end{displaymath}
\begin{displaymath}
\cong \left(\bigoplus_{\lambda\vdash n} d_{\lambda}  (V^{\otimes n}\otimes_{KS_n} S^{\lambda})\right) \oplus \left( \bigoplus_{(\lambda_1, \lambda_2)\vdash n} m(\lambda_1, \lambda_2) ( V^{\otimes n}\otimes_{KS_n} S^{(\lambda_1, \lambda_2)})\right)\cong
\end{displaymath}
\begin{displaymath}
\cong \left(\bigoplus_{\lambda\vdash n} d_{\lambda}  W^{\lambda}\right) \oplus \left( \bigoplus_{(\lambda_1, \lambda_2)\vdash n} m(\lambda_1, \lambda_2) W^{(\lambda_1, \lambda_2)}\right).
\end{displaymath}

{\bf d.} ($A_n$-case) As in case {\bf c} ( $S_n$-case )

{\bf e.} Follows from {\bf a} and Corollary 7.13.9 in ~\cite{Stan}  \ \ \ \ $\Box$

\section*{Acknowledgements}

$\bullet$ This project partially was carried out when the second named author visited the Institute of Mathematics and Informatics of the Bulgarian Academy of Sciences. He is very grateful for the creative atmosphere and the warm hospitality during his visit.

$\bullet$ Authors are grateful to V. Drensky for his kind interest in our results and for essential comments.

\vspace{1cm}



\end{document}